\documentclass{amsart}
\usepackage{amssymb,euscript,amsmath, mathrsfs}
\usepackage[dvips]{graphicx}
\usepackage[dvips]{color}

\newcounter{ENUM}
\newcommand{\itm}{\item}
\newenvironment{ilist}{\renewcommand{\theENUM}{\roman{ENUM}}\renewcommand{\itm}{\addtocounter{ENUM}{1}\item[(\theENUM)]}\begin{itemize}\setcounter{ENUM}{0}}{\end{itemize}}
\newenvironment{Ilist}{\renewcommand{\theENUM}{\Roman{ENUM}}\renewcommand{\itm}{\addtocounter{ENUM}{1}\item[(\theENUM)]}\begin{itemize}\setcounter{ENUM}{0}}{\end{itemize}}

\def\risom{\overset{\sim}{\rightarrow}}

\input xy
\xyoption{all}
\CompileMatrices

\def\GG{{\mathbb G}}

\def\cE{{\mathscr E}}

\def\cG{{\mathscr G}}

\def\cL{{\mathscr L}}
\def\cO{{\mathscr O}}

\def\cM{{\mathcal M}}

\def\cLG{\mathcal{LG}}

\def\ch{\operatorname{char}}

\def\Spec{\operatorname{Spec}}
\def\Pic{\operatorname{Pic}}

\def\im{\operatorname{im}}

\def\EH{\operatorname{EH}}

\def\pr{\operatorname{pr}}
\def\d{\operatorname{d}}

\newcommand{\margh}[1]{}

\newtheorem{thm}{Theorem}[section]

\newtheorem{lem}[thm]{Lemma}
\newtheorem{cor}[thm]{Corollary}

\theoremstyle{definition}
\newtheorem{defn}[thm]{Definition}
\newtheorem{ques}[thm]{Question}

\newtheorem{sit}[thm]{Situation}
\theoremstyle{remark}
\newtheorem{notn}[thm]{Notation}
\newtheorem{rem}[thm]{Remark}

\numberwithin{equation}{section}

\begin{document}
\title{Linked Grassmannians and crude limit linear series}
\author{Brian Osserman}
\begin{abstract} 
In \cite{os8}, a new construction of limit linear series is presented
which functorializes and compactifies the original construction of
Eisenbud and Harris, using a new space called the linked Grassmannian. The 
boundary of the compactification consists of crude limit series,
and maps with positive-dimensional fibers to crude limit series of
Eisenbud and Harris. In this paper, we carry out a careful analysis of
the linked Grassmannian to obtain an upper bound on the dimension of the
fibers of the map on crude limit series, thereby concluding an upper bound on 
the dimension of the locus of crude limit series, and obtaining a simple 
proof of the 
Brill-Noether theorem using only the limit linear series machinery. We also
see that on a general reducible curve, even crude limit series may be 
smoothed to nearby fibers.
\end{abstract}
\thanks{The author was supported by a fellowship from the National Science 
Foundation during the preparation of this paper.}
\maketitle

\section{Introduction}

In \cite{e-h1}, Eisenbud and Harris introduced the powerful theory of 
limit linear series. In \cite{os8}, a new construction of spaces 
of limit linear series is introduced,
which is functorial and provides a compactification of the Eisenbud-Harris
space, agreeing with it on the open locus of ``refined'' limit series. 
Coincidentally, concrete motivation for such a construction has been 
provided by Khosla \cite{kh1}, who produces, given a suitable proper stack 
of limit linear
series, an infinite family of effective virtual divisors in $\cM_g$, and
shows that whenever they are divisors, they give counter-examples to the
Harris-Morrison slope conjecture.

The Eisenbud-Harris construction is forced to omit what they call the locus
of ``crude'' limit series; they give a fiber-by-fiber description of this
boundary, but do not include it in the relative construction which is the
heart of the theory. The theory of \cite{os8} has a boundary which maps
naturally to the Eisenbud-Harris crude limit series, but frequently with
positive-dimensional fibers. Because of this distinction, we will refer to
the boundary elements of the latter construction as crude limit series,
and the boundary described by Eisenbud and Harris as $EH$-crude limit 
series. 

Dimension estimates are central to both theories of limit linear
series, and while $EH$-crude limit series are easily amenable to making such
estimates inductively, the crude limit series of \cite{os8} are more 
combinatorially
complicated, and were not closely analyzed in \cite{os8}. The goal of the
present paper is to obtain sufficiently sharp upper bounds for the 
dimensions of spaces of crude limit series that we can apply the theoretical
machinery of \cite{os8} to the loci of crude limit series in addition to
refined limit series. Our estimates will allow us to prove the following
theorem.

\begin{thm}\label{main} Fix integers $r,d$, and let $X$ be a general curve 
of compact type over $\Spec k$ having no more than two components, with 
$\ch k = 0$, and general marked points. Then the space of limit linear series 
on $X$ of degree $d$ and dimension $r$, with prescribed ramification at the 
marked points, is proper, and pure of exactly the expected dimension 
$\rho=(r+1)(d-r)-rg-\sum_{i,j}\alpha^i_j$.

If no ramification is specified, this space is non-empty, and if further
$\rho>0$, the space is connected.
\end{thm}

Here by a general curve of compact type, we mean that the dual graph and
the genus of each component may be specified, and then the isomorphism class
of each component must be allowed to be general. In fact, everything except 
the connectedness will follow easily from the theory of 
\cite{os8}, so in particular we obtain a new and direct proof of the
Brill-Noether theorem for linear series with prescribed ramification
\cite[Thm. 4.5]{e-h1} in characteristic $0$. We also see that on a general
curve, all limit series, including crude limit series, are smoothable to 
nearby fibers. We thus obtain positive answers to Questions 7.1 and 7.3 of
\cite{os8}. 

Unfortunately, the presence of 
inseparable linear series poses an obstacle to carrying through the same
proof in positive characteristic, although a different proof for the
case of $1$-dimensional linear series is given in 
\cite{os3}. Finally, we mention the complementary result 
\cite[Thm. 4.3]{o-h1} that when a limit linear space has the expected 
dimension, it is Cohen-Macaulay and flat over the base. We therefore
conclude that over a general curve with two components in characteristic 0, 
limit linear series schemes are quite well behaved.

See \S \ref{s-apps} below for background on limit linear series and
smoothing families. The main theorem
is proved by careful analysis of the linked Grassmannian, which arises
in the limit linear series construction of \cite{os8}. 
We begin in \S \ref{s-review} by reviewing the definition of and basic 
results on the linked Grassmannian. We then focus our attention in \S
\ref{s-back} on the map which projects to the first and last subspaces, 
introducing some notation and stating background lemmas to set up a more
detailed analysis of the fibers of this map, which is carried out in
\S \ref{s-fibers}. Finally, we carry out the stated application to
spaces of limit linear series in \S \ref{s-apps}.

\section{Review of the linked Grassmannian}\label{s-review}

We briefly review the basic definitions and results of the linked 
Grassmannian. 

\begin{defn}\label{lgdef} Let $S$ be an integral, locally Cohen-Macaulay 
scheme, and $\cE_1, \dots, \cE_n$ vector 
bundles on $S$, each of rank $d$. Given maps
$f_i:\cE_i \rightarrow \cE_{i+1}$ and $g_i: \cE_{i+1} \rightarrow
\cE_i$, and a positive integer $r<d$, we denote by $\cLG:=\cLG(r,\{\cE_i\}_i,
\{f_i, g_i\}_i)$ the functor associating to each $S$-scheme $T$ the set of 
sub-bundles $V_1, \dots, V_n$ of $\cE_{1,T}, \dots, \cE_{n,T}$ having rank 
$r$ and 
satisfying $f_{i,T}(V_i) \subset V_{i+1}$, $g_{i,T}(V_{i+1}) \subset V_i$ 
for all $i$.

We say that $\cLG$ is a {\bf linked Grassmannian} functor if the following
further conditions on the $f_i$ and $g_i$ are satisfied:

\begin{Ilist}
\itm There exists some $s \in \cO_S$ such that $f_i g_i = g_i f_i$ is
scalar multiplication by $s$ for all $i$.
\itm Wherever $s$ vanishes, the kernel of $f_i$ is precisely equal to 
the image of $g_i$, and vice versa. More precisely, for any $i$ and given any
two integers $r_1$ and $r_2$ such that $r_1 + r_2 < d$, then the closed
subscheme of $S$ obtained as the locus where $f_i$ has rank less than or
equal to $r_1$ and $g_i$ has rank less than or equal to $r_2$ is empty.
\itm At any point of $S$, $\im f_i \cap \ker f_{i+1}=0$, and $\im
g_{i+1} \cap \ker g_i = 0$. More precisely, for any integer $r_1$, and any
$i$, we have locally closed subschemes of $S$ corresponding to the locus
where $f_i$ has rank exactly $r_1$, and $f_{i+1} f_i$ has rank less than or
equal to $r_1-1$, and similarly for the $g_i$. Then we require simply that
all of these subschemes be empty.
\end{Ilist}
\end{defn}

The main theorem of \cite{os8} on the linked Grassmannian is the following:

\begin{thm}\label{lg-main} \cite[Lem. A.3, Thm. A.15]{os8} $\cLG$ is 
representable by a 
scheme $LG$; this scheme is naturally a closed subscheme of the obvious 
product $G_1 \times \cdots \times G_n$ of Grassmannian schemes over $S$, 
which is smooth of relative dimension $nr(d-r)$. Each component of $LG$ has 
codimension $(n-1)r(d-r)$ in the product, and maps
surjectively to $S$. If $s$ is non-zero, then $LG$ is also irreducible.
\end{thm}

We consider the following question, motivated by applications to the theory
of limit linear series:

\begin{ques} Let $LG=LG(E_i,f_i,g_i)$ be a linked Grassmannian space over 
$\Spec k$, and $\pr_{1n}: LG \rightarrow G_1 \times G_n$ the projection map 
determined by forgetting all but the first and last subspaces. What are the 
dimensions of the fibers of this map?
\end{ques}

We thus fix $(V_1,V_n) \in G_1 \times G_n$, and want to consider the ways
of filling in the intermediate $V_i\subset E_i$ to obtain a collection 
of spaces linked by the $f_i$ and $g_i$. Because we are motivated primarily
by applications to limit linear series, we will in fact be interested 
primarily in obtaining upper bounds rather than computing the precise 
answer. Furthermore, because the problem is otherwise trivial, we fix
the following assumptions:

\begin{sit} We have fixed vector spaces $E_i$ and maps $f_i:E_i
\rightarrow E_{i+1}$, and $g_i:E_{i+1}\rightarrow E_i$ satisfying
the hypotheses of Definition \ref{lgdef}. We further assume that $n>2$ 
and the $s$ of condition (I) Definition \ref{lgdef} is equal
to $0$. Finally, we fix $V_1 \subset E_1$ and $V_n \subset E_n$ such that
the iterated image of $V_1$ under the $f_i$ is contained in $V_n$, and
the iterated image of the $V_n$ under the $g_i$ is contained in $V_1$.
\end{sit}

We remark that although the last condition is certainly necessary for the
pair $(V_1,V_n)$ to come from a point of the linked Grassmannian, it is by
no means sufficient. In the course of our analysis we will produce a 
sufficient condition, albeit an extremely unwieldy one.

\section{Notation and background lemmas}\label{s-back}

We begin with some preliminary definitions and observations.

\begin{notn}\label{nfirst} For $1 \leq i \leq j \leq n-1$, we denote by 
$f_{i,j}$ the 
composition $f_j \circ \dots \circ f_i$, and $g_{j,i}$ the composition 
$g_i \circ \dots \circ g_j$. 

Given $V_1 \subset E_1$ and $V_n \subset E_n$, 
for $1 \leq i \leq n$, we write:
$$\bar{V}_{1,i}:= g_{i-1,1}(g_{i-1,1}^{-1}(V_1))=V_1 \cap \im g_{i-1,1}
\subset V_1, \text{ and}$$
$$\bar{V}_{n,i} := f_{i,n-1}(f_{i,n-1}^{-1}(V_1))=V_1 \cap \im f_{i,n-1}
\subset V_n,$$
where by convention $g_{0,1}$ and $f_{n,n-1}$ are just the identity map,
so that $\bar{V}_{1,1}=V_1$, and $\bar{V}_{n,n}=V_n$.
We also write
$$V_{1,n}:=g_{n-1,1}(V_n),\text{ and }
V_{n,1}:= f_{1,n-1}(V_1).$$
\end{notn}

We thus obtain filtrations 
$$V_{1,n} \subset \bar{V}_{1,n} \subset \bar{V}_{1,n-1}
\subset \dots \subset \bar{V}_{1,2} \subset V_1$$ 
and
$$V_{n,1} \subset \bar{V}_{n,1} \subset \bar{V}_{n,2} 
\subset \dots \subset \bar{V}_{n,n-1} \subset V_n.$$ 
Note that the first containments of each filtration make use
of our hypothesis that $g_{n-1,1}(V_n)\subset V_1$ and $f_{1,n-1}(V_1)\subset
V_n$. 

We now observe that any intermediate set of 
$V_i$ linking $V_1$ and $V_n$ can be constructed within the spaces 
$\bar{V}_{1,i}\oplus \bar{V}_{n,i}$.

\begin{lem}\label{lem-1} Let $\{V_i\}$ be a point of a linked Grassmannian 
$LG(E_i,f_i,g_i)$. Then for each $i=2,\dots,n-1$ there is a natural 
injection
$$V_i \hookrightarrow g_{i-1,1}^{-1}(V_1) \cap f_{i,n-1}^{-1}(V_n) 
\hookrightarrow \bar{V}_{1,i} \oplus \bar{V}_{n,i},$$ 
defined by the map $g_{i-1,1}\oplus f_{i,n-1}$.
\end{lem}

\begin{proof} The first inclusion inside $E_i$ follows from the assumption
that the $V_i$ are all linked under the $f_i$ and $g_i$. For the second,
we only need to see that any vector in $E_i$ mapping to $0$ under 
$g_{i-1,1}$ and $f_{i,n-1}$ must itself be $0$. By condition (III) of
Definition \ref{lgdef}, we have $\ker g_{i-1,1} = \ker g_{i-1}$ and 
$\ker f_{i,n-1} = \ker f_i$, and further $\ker f_i \cap \im f_{i-1} = 0$.
But by condition (II) of {\it loc.\ cit.}, $\ker g_{i-1} = \im f_{i-1}$
and is hence disjoint from $\ker f_i$, completing the proof.
\end{proof}

\begin{notn} In the situation of the lemma, we denote by $\bar{Z}_i$ the 
cokernel of the map $g_{i-1,1}^{-1}(V_1) \cap f_{i,n-1}^{-1}(V_n) 
\hookrightarrow \bar{V}_{1,i} \oplus \bar{V}_{n,i}$, and by $\bar{Z}_{1,i}$ 
and $\bar{Z}_{n,i}$ the images of $\bar{V}_{1,i}\oplus (0)$ and $(0) \oplus
\bar{V}_{n,i}$, respectively.
\end{notn}

We next observe the following:

\begin{lem}\label{lem-2} For all $i=2,\dots,n-1$, we have
$$\bar{V}_{1,i+1}= \bar{V}_{1,i} \cap g_{i,1}(E_{i+1}),$$
$$\bar{V}_{n,i-1}= \bar{V}_{n,i} \cap f_{i-1,n-1}(E_{i-1}),$$
and short exact sequences
$$0 \rightarrow \bar{V}_{1,i+1} \rightarrow \bar{V}_{1,i} \rightarrow 
\bar{Z}_{1,i} \rightarrow 0,$$ 
$$0 \rightarrow \bar{V}_{n,i-1} \rightarrow \bar{V}_{n,i} \rightarrow 
\bar{Z}_{n,i} \rightarrow 0.$$ 
\end{lem}

\begin{proof} The first two equalities follow immediately from 
the definitions, as do injectivity and surjectivity of the sequences. For 
exactness of the first sequence, we check the following identities:
\begin{align*}\bar{V}_{1,i+1}& = \bar{V}_{1,i} \cap g_{i,1}(E_{i+1}) \\
&=g_{i-1,1}(g_{i-1,1}^{-1}(V_1)) \cap g_{i-1,1}(g_i(E_{i+1})) \\ 
&=g_{i-1,1}(g_{i-1,1}^{-1}(V_1) \cap g_i(E_{i+1})) \\ 
&=g_{i-1,1}(g_{i-1,1}^{-1}(V_1) \cap \ker f_i) \\ 
&= \ker(\bar{V}_{1,i} \to \bar{Z}_{1,i}),
\end{align*}
and exactness of the second sequence follows similarly.
\end{proof}

\begin{notn} Given a point $\{V_i\}$ of $LG$, we denote by 
$V_{1,i}\subset V_1$ and $V_{n,i}\subset V_n$ the images $g_{i-1,1}(V_i)$ 
and $f_{i,n-1}(V_i)$ respectively, and by $Z_i \subset \bar{Z}_i$ 
the image of $V_{1,i} \oplus V_{n,i}$ in
$\bar{Z}_i$. 

Finally, for any vector space $V$, if we wish to prescribe
a certain dimension for $V$, we will fix a non-negative integer which we will
denote by $d_V$, and require that $V$ satisfy $\dim (V) = d_V$; we will
also abbreviate $\dim (V)$ by $\d(V)$. 
\end{notn}

Note that this is compatible with our earlier definitions of $V_{1,n}$,
$V_{n,1}$. 

\begin{lem}\label{lem-3} With $V_1, V_n$ given, let $V_{1,n}$ and $V_{n,1}$
be as determined by $(V_1,V_n)$, and fix nested sequences of subspaces 
$$V_{1,n} \subset V_{1,n-1} \subset ... \subset V_{1,2} \subset V_1,
\text{ and}$$
$$V_{n,1} \subset V_{n,2} \subset ... \subset V_{n,n-1} \subset V_n$$
with $V_{1,i}\subset \bar{V}_{1,i}$ and $V_{n,i} \subset \bar{V}_{n,i}$
for all $i$. Then points of the linked Grassmannian with the given $V_1,V_n$,
and $V_{1,i}$ and $V_{n,i}$ correspond to collections of $r$-dimensional 
subspaces $V_i \subset \ker (V_{1,i}\oplus V_{n,i} \rightarrow \bar{Z}_i)$ 
with $V_i$ mapping surjectively to $V_{1,i}$ and $V_{n,i}$ and containing
$V_{1,i+1}\oplus V_{n,i-1}$ for $i=2,\dots,n-1$.
\end{lem}

\begin{proof} By Lemma \ref{lem-1}, any $\{V_i\}$ corresponding to a point 
of the
linked Grassmannian with given associated $V_{1,i}$ and $V_{n,i}$ may be 
considered as subspaces of
$\ker (V_{1,i}\oplus V_{n,i} \rightarrow \bar{Z}_i)$, surjecting on $V_{1,i}$
and $V_{n,i}$. Conversely any $V_i$ of this form may naturally be considered 
as a subspace of $E_i$. Thus, we need only check 
that the linking condition is equivalent to
the containment of $V_{1,i+1} \oplus V_{n,i-1}$, which by Lemma \ref{lem-2}
is contained in the kernel space in question.

We show that linking under the $f_i$ is equivalent to $V_i$ containing
$(0) \oplus V_{n,i-1}$, with the argument for the $g_i$ being the same, and
the combination giving the desired statement. First, 
$f_1 (V_1) \subset V_2$: indeed, the linked Grassmannian conditions and
our definitions imply that 
$f_1(V_1) \risom (0) \oplus V_{n,1} \subset E_1 \oplus E_n$, and 
$(0) \oplus V_{n,1} \subset V_1 \oplus V_n$ by hypothesis. Hence, 
$f_1(V_1)$ is contained in
$V_2$ if and only if $V_2$ contains $(0) \oplus V_{n,1}$.
Similarly, for $2<i\leq n-1$, we have $f_{i-1}(V_{i-1}) \subset V_i$
if and only if $V_i$ contains $(0) \oplus V_{n,i-1}$. Lastly, the
containment $f_{n-1}(V_{n-1})=V_{n,n-1}\subset V_n$ is automatic from our 
hypotheses.
\end{proof}

\begin{lem}\label{lem-4} Given a point $\{V_i\}$ of $LG$, each
$V_{1,i}$ and $V_{n,i}$ maps surjectively to $Z_i$, so in particular 
$Z_i \subset \bar{Z}_{1,i} \cap \bar{Z}_{n,i}$.
\end{lem}

\begin{proof} This is equivalent to the statement that the images of 
$V_{1,i}$ and $V_{n,i}$ agree inside $\bar{Z}_i$. But for any 
$v \in V_{1,i}$, by definition there is a vector 
$(v,v')\in V_{1,i} \oplus V_{n,i} \subset \bar{V}_{1,i} \oplus \bar{V}_{n,i}$ 
coming from $V_i$, hence mapping to $0$ in $Z_i$. Thus $-v'$ has the same 
image as $v$ in $\bar{Z}_i$, so the image of $V_{1,i}$ is contained in that
of $V_{n,i}$. But the same argument works to show the opposite containment,
giving the desired statement.
\end{proof}

The final technical lemma is the following:

\begin{lem}\label{strata} Fixing non-negative integer dimensions 
$d_{V_{1,i}}, d_{V_{n,i}}, d_{Z_i}$ for all $i$ with $1<i<n$ determines
a locally closed subscheme of the fiber of the $\pr_{1n}$ map over 
$(V_1,V_n)$.
\end{lem}

\begin{proof} First, the conditions obtained by fixing the $d_{V_{1,i}}$
and $d_{V_{n,i}}$ determine a locally closed subscheme, since they are
imposing a particular rank on the maps of the universal bundles 
$V_i \overset{g_{i-1,1}}{\rightarrow} V_1$ and 
$V_i \overset{f_{i,n-1}}{\rightarrow} V_n$, which are locally free. Next,
within the locally closed subscheme cut out by these conditions, we note
that the $V_{1,i}$ and $V_{n,i}$ are also locally free (and the $\bar{Z}_i$, 
being determined by $(V_1,V_n)$, are in fact free), so prescribing
the ranks of the maps $V_{1,i} \rightarrow \bar{Z}_i$ (or equivalently, by 
Lemma \ref{lem-4}, the ranks of the maps $V_{n,i} \rightarrow \bar{Z}_i$) 
determines locally closed conditions.
\end{proof}

We thus obtain a stratification of the fiber of $\pr_{1n}$, in the sense
of having a collection of disjoint locally closed subschemes whose union
is set-theoretically the entire fiber. Our main task will be to analyze
this stratification further.

\section{Dimensions of the fibers}\label{s-fibers}

We will initially analyze the pieces of our stratification to compute their
dimensions. Altough we describe the fibers of $\pr_{1n}$ in terms of the 
dimensional invariants 
$\d(\bar{V}_{1,i}),\d(\bar{V}_{n,i}),\d(\bar{Z}_i)$, 
the formulas are quite complicated, and rather than work with them directly,
we will be able to obtain an indirect bound by studying the dimension of 
pairs $(V_1,V_n)$ having specified dimensional invariants. We can then use 
the fact that the dimensions of the fibers are determined entirely by
the numerical invariants to obtain an indirect upper bound for them,
Corollary \ref{fiber-dim} below. 

\begin{thm}\label{stratum-dim} The fibers of $\pr_{1n}$ have dimension 
determined by the
dimensions $\d(V_{1,n})$, $\d(V_{n,1})$, $\d(\bar{V}_{1,n})$, 
$\d(\bar{V}_{n,1})$, and 
$\d(\bar{V}_{1,i})$, $\d(\bar{V}_{n,i})$, and $\d(\bar{Z}_i)$ for 
$1<i<n$.

Specifically, assuming $(V_1,V_n)$ satisfies $f_{1,n-1}(V_1)\subset V_n,
g_{n-1,1}V_n \subset V_1$, if we prescribe dimensions $d_{V_{1,i}}, 
d_{V_{n,i}}, d_{Z_i}$, the corresponding stratum of Lemma \ref{strata} is
non-empty precisely when the following conditions are satisfied for all $i$
with $2 \leq i \leq n-1$:

\begin{gather}
\label{cond-1} 
d_{Z_i}\leq \d(\bar{Z}_{1,i} \cap \bar{Z}_{n,i}) \\
\label{cond-2} 
\d(\bar{V}_{1,i+1}) \geq d_{V_{1,i}}-d_{Z_i} \geq d_{V_{1,i+1}}\\
\label{cond-3} 
\d(\bar{V}_{n,i-1}) \geq d_{V_{n,i}}-d_{Z_i} \geq d_{V_{n,i-1}} \\
\label{cond-4} 
d_{V_{1,i}}+d_{V_{n,i}}-d_{Z_i} \geq r \\
\label{cond-5} 
r \geq d_{V_{1,i+1}}+d_{V_{n,i}} \\
\label{cond-6} 
r \geq d_{V_{1,i}}+d_{V_{n,i-1}} 
\end{gather}
where by convention we write $d_{V_{n,1}}:=d(V_{n,1})$ and 
$d_{V_{1,n}}:=d(V_{1,1})$.

Each stratum is then smooth of dimension
\begin{multline}
\sum_{i=2}^{n-1} d_{Z_i} (\d(\bar{Z}_{1,i} \cap \bar{Z}_{n,i})-d_{Z_i}) \\
+ (d_{V_{1,i}}-d_{V_{1,i+1}})(\d(\bar{V}_{1,i+1})-d_{V_{1,i}}+d_{Z_i})
+ (d_{V_{n,i}}-d_{V_{n,i-1}})(\d(\bar{V}_{n,i-1})-d_{V_{n,i}}+d_{Z_i}) \\
+ (r-d_{V_{1,i+1}}-d_{V_{n,i-1}})(d_{V_{1,i}}+d_{V_{n,i}}-d_{Z_i}-r).
\end{multline}
\end{thm}

\begin{proof} We simply have to check that the dimensions (and emptyness
or non-emptyness) of the pieces of our stratification are determined by
the numbers in question. Note that by construction, for $2 \leq i \leq n-1$,
$\bar{Z}_{1,i}+\bar{Z}_{n,i}=\bar{Z}_i$, so 
$$\d(\bar{Z}_{1,i}\cap \bar{Z}_{n,i})= \d(\bar{Z}_{1,i}) + 
\d(\bar{Z}_{n,i}) - \d(\bar{Z}_i).$$ 
Furthermore, by Lemma \ref{lem-2} we have 
$$\d(\bar{Z}_{1,i})=\d(\bar{V}_{1,i})-\d(\bar{V}_{1,i+1}), \text{ and }
\d(\bar{Z}_{n,i})= \d(\bar{V}_{n,i})-\d(\bar{V}_{n,i-1}).$$ 
We fix non-negative integers $d_{V_{1,i}}, d_{V_{n,i}}, d_{Z_i}$, and
consider the structure of the resulting stratum of Lemma \ref{strata},
which we denote by $X$. 

Denote by $X'$ the functor of filtrations
$$g_{n-1,1}(V_n)=V_{1,n} \subset V_{1,n-1} \subset ... \subset V_{1,2} 
\subset V_1\text{ and}$$
$$f_{1,n-1}(V_1)=V_{n,1} \subset V_{n,2} \subset ... \subset V_{n,n-1} 
\subset V_n$$
with each $V_{1,i}\subset \bar{V}_{1,i}$ and $V_{n,i}\subset \bar{V}_{n,i}$,
and each of the prescribed dimension, with $V_{1,i}$ and $V_{n,i}$ having
the same image, also of the prescribed dimension, in $\bar{Z}_i$.
Also denote by $X''$ the functor of $(n-2)$-tuples of spaces 
$Z_i\subset \bar{Z}_{1,i}\cap \bar{Z}_{n,i}$ having the prescribed dimension. 
Then we have maps $X \rightarrow X' \rightarrow X''$,
which we will analyze one by one.

$X''$ is simply a product of Grassmannians, hence a smooth scheme 
of dimension 
$\sum_{i=2}^{n-1} d_{Z_i} (\d(\bar{Z}_{1,i} \cap \bar{Z}_{n,i})-d_{Z_i})$,
non-empty if and only if (\ref{cond-1}) is satisfied. 

We next wish to show that $X'$ is likewise representable by a scheme smooth 
over $X''$, of relative dimension
\begin{multline*} \sum_{i=n-1}^2 (d_{V_{1,i}}-d_{V_{1,i+1}})
(\d(\bar{V}_{1,i})-d_{V_{1,i}} -\d(\bar{Z}_{1,i})+d_{Z_i}) \\
+ \sum_{i=2}^{n-1} (d_{V_{n,i}}-d_{V_{n,i-1}})(\d(\bar{V}_{n,i})-d_{V_{n,i}}
-\d(\bar{Z}_{n,i})+d_{Z_i}),
\end{multline*} 
and non-empty if and only if (\ref{cond-2}) and (\ref{cond-3}) are satisfied
for $i$ with $2 \leq i \leq n-1$. Having chosen the $Z_i$,
the choices of $V_{1,i}$ are independent of the choices of 
$V_{n,i}$, and the situation is symmetric for both, so we treat only the 
former. We can describe $X'$ as open inside
a sequence of Grassmannian bundles over $X''$, with the first bundle 
corresponding to the choice of $V_{1,n-1}$, and each subsequent one
corresponding to the choice of some $V_{1,i}$ given $V_{1,i+1}$. 

Indeed, $V_{1,n-1}$ may be any space containing $V_{1,n}$, contained 
in the subspace of $\bar{V}_{1,n-1}$ mapping into $Z_{n-1}$, and surjecting 
onto $Z_{n-1}$. The last is an open condition, and we claim we have 
non-emptiness of dimension 
\begin{multline*} 
(d_{V_{1,n-1}}-\d(V_{1,n}))(\d(\bar{V}_{1,n-1})-d_{V_{1,n-1}}
-(\d(\bar{Z}_{1,n-1})-d_{Z_{n-1}}))\\
= (d_{V_{1,n-1}}-\d(V_{1,n}))(\d(\bar{V}_{1,n})-d_{V_{1,n-1}}
+d_{Z_{n-1}})
\end{multline*}
exactly when (\ref{cond-2}) is satisfied for $i=n-1$. Non-emptyness of the
Grassmannian is equivalent to $\d(\bar{V}_{1,n}) +d_{Z_{n-1}} \geq 
d_{V_{1,n-1}} \geq d_{V_{1,n}}$. We note that by Lemma \ref{lem-2}
we have that $V_{1,n}$ maps to $0$ in $\bar{Z}_i$, so
non-emptyness of the open condition of surjecting onto $Z_{n-1}$ is
then equivalent to $d_{V_{1,n-1}}- d_{V_{1,n}}\geq d_{Z_i}$, completing the
proof of the claim.
Similarly,
for each $i<n-1$, we can choose $V_{i,1}$ to be any space containing
$V_{1,i+1}$, contained in the subspace of $\bar{V}_{1,i}$ mapping into
$Z_i$, and surjecting onto $Z_i$. This is (open inside) another 
Grassmannian bundle, of dimension 
$(d_{V_{1,i}}-d_{V_{1,i+1}})(\d(\bar{V}_{1,i})-d_{V_{1,i}}
-\d(\bar{Z}_{1,i})+d_{Z_{1,i}})$
and non-empty exactly when (\ref{cond-2}) is satisfied for $i=n-2,\dots,2$.
Arguing the same way for the $V_{n,i}$, we find that (\ref{cond-3})
is precisely the condition for non-emptyness, and we conclude the desired 
description of $X'$.

Finally, given the $V_{1,i}$ and $V_{n,i}$, by Lemma \ref{lem-3} we have
that $X$ is open inside a bundle of products of Grassmannians of dimension 
$$\sum_{i=2}^{n-1} (r-d_{V_{1,i+1}}-d_{V_{n,i-1}})
(d_{V_{1,i}}+d_{V_{n,i}}-d_{Z_i}-r)$$
over $X'$. Now, non-emptyness
of the bundle is equivalent to (\ref{cond-4}) together with
$$r \geq d_{V_{1,i+1}}+d_{V_{n,i-1}},$$
and we claim that non-emptyness of the
open condition that each $V_i$ surject onto $V_{1,i}$ and $V_{n,i}$
is equivalent to (\ref{cond-5}) and (\ref{cond-6}). 
This claim is easily checked, recalling the conditions imposed on $V_i$ by
the containment of $V_{1,i+1} \oplus (0)$ and $(0) \oplus V_{n,i-1}$, by the 
observation that 
because each of $V_{1,i}$ and $V_{n,i}$ surject onto $Z_i$, we must have that
$\ker (V_{1,i}\oplus V_{n,i} \rightarrow Z_i)$ surjects onto $V_{1,i}$ and
$V_{n,i}$. Finally, either of (\ref{cond-5}) or (\ref{cond-6}) implies the 
above
inequality, so we find that non-emptiness of $X$, given non-emptyness of $X'$,
is equivalent to (\ref{cond-4}), (\ref{cond-5}), and (\ref{cond-6}).

This shows that the dimension and non-emptyness of $X$ are 
entirely determined by the various prescribed dimensions, and are otherwise 
independent of $(V_1,V_n)$, completing the proof of the theorem.
\end{proof}

Because of the complexity of the formula obtained from the preceding theorem,
we will find it profitable to analyze the situation less directly, by
considering spaces of spaces $(V_1,V_n)$ with given dimensional invariants.
We first need to introduce one more piece of 
notation:

\begin{notn} We denote by $\tilde{Z}_i$ the cokernel of the map
$$E_i \hookrightarrow f_{i,n-1}(E_i) \oplus g_{i-1,1}(E_i).$$
\end{notn}

Note that the injectivity of the map follows from the same argument as in
Lemma \ref{lem-1}. Note also that directly from the definition, one sees 
that each
of $f_{i,n-1}(E_i)$ and $g_{i-1,1}(E_i)$ surjects onto $\tilde{Z}_i$.

\begin{thm}\label{crude-loci} The dimension of the space of pairs 
$(V_1,V_n) \in G_1 \times G_n$
having dimensions prescribed by fixed 
$d_{\bar{V}_{1,i}}, d_{\bar{V}_{n,i}}$ and 
$d_{\bar{Z}_i}$ for $i=2,\dots,n-1$, as well as 
$d_{V_{1,n}}, d_{V_{n,1}}, d_{\bar{V}_{1,n}}$, and $d_{\bar{V}_{n,1}}$,
and satisfying 
$f_{1,n-1}(V_1) \subset V_n$, $g_{n-1,1}(V_n)
\subset V_1$, is given by 
\begin{multline*}\sum_{i=1}^{n} d_{\bar{Z}_{1,i}}
(d_{g_{i-1,1}(E_i)}-d_{\bar{V}_{1,i}})+
d_{\bar{Z}_{n,i}} (d_{f_{i,n-1}(E_i)}-d_{\bar{V}_{n,i}}) \\
+d_{V_{1,n}}(d_{\bar{V}_{1,n}}-d_{V_{1,n}})
+d_{V_{n,1}}(d_{\bar{V}_{n,1}}-d_{V_{n,1}})
- \sum_{i=2}^{n-1} d_{\bar{Z}_{1,i}\cap\bar{Z}_{n,i}}
(d_{\tilde{Z}_i}-d_{\bar{Z}_i}),
\end{multline*}
where we set
$$d_{\bar{Z}_{1,i}}=d_{\bar{V}_{1,i}}-d_{\bar{V}_{1,i+1}}, \text{ and }
d_{\bar{Z}_{n,i}}= d_{\bar{V}_{n,i}}-d_{\bar{V}_{n,i-1}},$$ 
and
$$d_{\bar{Z}_{1,i}\cap \bar{Z}_{n,i}}= d_{\bar{Z}_{1,i}} + 
d_{\bar{Z}_{n,i}} - d_{\bar{Z}_i}$$ 
for each $i=2,\dots,n-1$,
and use the conventions that 
$$d_{\bar{V}_{1,n+1}}=d_{\bar{V}_{n,0}}=0,$$
$$d_{\bar{V}_{1,1}}=d_{\bar{V}_{n,n}}=r,$$ 
$$d_{g_{0,1}(E_1)}=d_{V_{n,1}}+d_{\ker f_1}, \text{ and}$$ 
$$d_{f_{n,n-1}(E_1)}=d_{V_{1,n}}+d_{\ker g_{n-1}}.$$ 

This space is smooth, and non-empty if and only if the following conditions 
are satisfied for $i=2,\dots,n-1$:

\begin{gather}
\label{2-cond-1} 
d_{\bar{Z}_{1,i}\cap\bar{Z}_{n,i}} \leq d_{\bar{Z}_{1,i}}
\leq d_{\tilde{Z}_i} \\
\label{2-cond-2} 
d_{\bar{Z}_{1,i}\cap\bar{Z}_{n,i}} \leq d_{\bar{Z}_{n,i}}
\leq d_{\tilde{Z}_i} \\
\label{2-cond-3}
d_{\bar{Z}_{1,i}}+ d_{\bar{Z}_{n,i}}
\leq d_{\tilde{Z}_i} + d_{\bar{Z}_{1,i}\cap\bar{Z}_{n,i}} \\
\label{2-cond-4}
d_{V_{n,1}} \leq d_{\bar{V}_{n,1}} \leq d_{f_{1,n-1}(E_1)} \\
\label{2-cond-5}
d_{V_{1,n}} \leq d_{\bar{V}_{1,n}} \leq d_{g_{n-1,1}(E_n)} \\
\label{2-cond-6}
d_{\bar{V}_{1,i+1}} \leq d_{\bar{V}_{1,i}} \leq 
d_{g_{i-1,1}(E_i)}+d_{\bar{Z}_{1,i}}-d_{\tilde{Z}_i} \\
\label{2-cond-7}
d_{\bar{V}_{n,i-1}} \leq d_{\bar{V}_{n,i}} \leq 
d_{f_{i,n-1}(E_i)}+d_{\bar{Z}_{n,i}}-d_{\tilde{Z}_i} \\
\label{2-cond-8}
d_{\bar{V}_{1,2}}+d_{V_{n,1}} \leq d_{\bar{V}_{1,1}} \leq 
d_{g_{0,1}(E_1)} \\
\label{2-cond-9}
d_{\bar{V}_{n,n-1}}+d_{V_{1,n}} \leq d_{\bar{V}_{n,n}} \leq 
d_{f_{n,n-1}(E_n)} 
\end{gather}
\end{thm} 

\begin{proof} The proof proceeds similarly to that of Theorem 
\ref{stratum-dim}: we build up $V_1$ by one $\bar{V}_{1,i}$ at a time,
and similarly for $V_n$. However, we will have to begin by 
choosing the $\bar{Z}_{1,i}\cap \bar{Z}_{n,i}$, followed by the 
$\bar{Z}_{1,i}$ and
$\bar{Z}_{n,i}$. If we denote by $X$ the space of pairs $(V_1,V_n)$ of the
appropriate form, we will want to consider the functors $X'$ of 
$(2n-4)$-tuples of $\bar{Z}_{1,i}$ and $\bar{Z}_{n,i}$ inside $\tilde{Z}_i$
such that $\d (\bar{Z}_{1,i} \cap \bar{Z}_{n,i}) = 
d_{\bar{Z}_{1,i}\cap \bar{Z}_{n,i}}$, and $X''$ of $(n-2)$-tuples of
$\bar{Z}_{1,i}\cap\bar{Z}_{n,i}$ inside $\tilde{Z}_i$ of dimension
$d_{\bar{Z}_{1,i}\cap \bar{Z}_{n,i}}$. As before, we have natural maps
$X \rightarrow X' \rightarrow X''$ which we analyze one at a time.

Now, $X''$ is clearly represented by a product of Grassmannians of total
dimension
$$\sum_{i=2}^{n-1} d_{\bar{Z}_{1,i}\cap \bar{Z}_{n,i}} (d_{\tilde{Z}_i}
- d_{\bar{Z}_{1,i}\cap \bar{Z}_{n,i}}),$$
non-empty if and only if 
$d_{\bar{Z}_{1,i}\cap\bar{Z}_{n,i}} \leq d_{\tilde{Z}_i}$, which is implied
by (\ref{2-cond-1}).
$X'$ is then open inside a bundle in products of Grassmannians over $X''$, 
of dimension 
$$\sum_{i=2}^{n-1} (d_{\bar{Z}_{1,i}}-d_{\bar{Z}_{1,i}\cap\bar{Z}_{n,i}})
(d_{\tilde{Z}_i}-d_{\bar{Z}_{1,i}}) +
(d_{\bar{Z}_{n,i}}-d_{\bar{Z}_{1,i}\cap\bar{Z}_{n,i}})
(d_{\tilde{Z}_i}-d_{\bar{Z}_{n,i}});$$
the bundle is non-empty if and only if (\ref{2-cond-1}) and (\ref{2-cond-2})
are satisfied, while the open condition that the intersection of
$\bar{Z}_{1,i}$ and $\bar{Z}_{n,i}$ is no larger than the space chosen for 
$\bar{Z}_{1,i}\cap \bar{Z}_{n,i}$ is non-empty if and only if
(\ref{2-cond-3}) is satisfied.
 
Finally, we show that $X$ is made up of a tower of open subschemes of 
Grassmannian bundles over $X'$, by building up $V_1$ and $V_n$ starting from 
$V_{1,n}$ and $V_{n,1}$ and continuing up through the 
$\bar{V}_{1,i}$ and $\bar{V}_{n,i}$. The conditions on the $\bar{V}_{1,i}$
and $\bar{V}_{n,i}$ as we build them up will be simply that for each 
$i=2,\dots, n-1$, we have 
$\bar{V}_{1,i} \cap g_{i,1}(E_{i+1})= \bar{V}_{1,i+1}$, and 
$\bar{V}_{n,i} \cap f_{i-1,n-1}(E_{i-1})= \bar{V}_{n,i-1}$. Of course, we
will also require that $V_{1,n}=g_{n-1,1}(V_n)$ and
$V_{n,1}=f_{1,n-1}(V_1)$. These conditions will 
ensure that the $\bar{V}_{1,i}$ and $\bar{V}_{n,i}$ in fact come
from $(V_1,V_n)$ as prescribed by Notation \ref{nfirst}.

The spaces $\bar{V}_{1,i}$ are almost independent from the $\bar{V}_{n,i}$, 
except for the final requirement on $V_{1,n}$ and $V_{n,1}$;
however, since all $\bar{V}_{1,i}$ for $i>1$ are by definition in the image 
of $g_1$, they map to $0$ under $f_1$ and hence $f_{1,n-1}$, and similarly
for the $\bar{V}_{n,i}$, so this dependence will only appear when we choose
$V_1=\bar{V}_{1,1}$ and $V_n=\bar{V}_{n,n}$,
after all the previous $\bar{V}_{1,i}$ and $\bar{V}_{n,i}$
have been chosen. 

First, choosing $V_{1,n}\subset g_{n-1,1}(E_n)$ and $V_{n,1}
\subset f_{1,n-1}(E_1)$ 
is clearly a product of Grassmannians of dimension 
$$d_{V_{1,n}}(d_{g_{n-1,1}(E_n)}-d_{V_{1,n}})+ 
d_{V_{n,1}}(d_{f_{1,n-1}(E_1)}-d_{V_{n,1}}),$$ 
non-empty as long as (\ref{2-cond-4}) and (\ref{2-cond-5}) are satisfied. 
Next, choosing $\bar{V}_{1,n}\subset g_{n-1,1}(E_n)$ containing
$V_{1,n}$, and $\bar{V}_{n,1} \subset f_{1,n-1}(E_1)$ containing $V_{n,1}$
is again a product of Grassmannians, of dimension 
$$(d_{\bar{V}_{1,n}}-d_{V_{1,n}})(d_{g_{n-1,1}(E_n)}-d_{\bar{V}_{1,n}})+ 
(d_{\bar{V}_{n,1}}-d_{V_{n,1}})(d_{f_{1,n-1}(E_1)}-d_{\bar{V}_{n,1}}),$$ 
non-empty as long as (\ref{2-cond-4}) and (\ref{2-cond-5}) are satisfied; 
moreover, we see that to have non-emptiness in both cases, we must have
(\ref{2-cond-4}) and (\ref{2-cond-5}).

For $i=n-1,\dots,2$, we allow $\bar{V}_{1,i}$ to be an arbitrary subspace
of the preimage of our chosen $\bar{Z}_{1,i}$ inside $g_{i-1,1}(E_i)$, which
must contain $\bar{V}_{1,i+1}$, map surjectively onto $\bar{Z}_{1,i}$, and
must intersect with $g_{i,1}(E_{i+1})$ in precisely $\bar{V}_{1,i+1}$.
Because each $g_{i-1,1}(E_i)$ surjects onto $\tilde{Z}_i$, 
this will be open inside a Grassmannian of dimension
$$(d_{\bar{V}_{1,i}}-d_{\bar{V}_{1,i+1}})(d_{g_{i-1,1}(E_i)}-
(d_{\tilde{Z}_i}-d_{\bar{Z}_{1,i}})-d_{\bar{V}_{1,i}}).$$ 
This Grassmannian is non-empty if and only if 
(\ref{2-cond-6}) is satisfied, so it remains to analyze the open conditions,
which we claim are always non-empty. 

Noting that $g_{i,1}(E_{i+1})=\ker(g_{i-1,1}(E_i)\to \tilde{Z}_i)$, 
one checks that the condition of surjecting onto $\bar{Z}_{1,i}$ is 
equivalent to the condition that $\bar{V}_{1,i}$ 
intersect $g_{i,1}(E_{i+1})$ exactly in $\bar{V}_{1,i+1}$, and the
non-emptiness of both is equivalent to the inequality
$d_{\bar{V}_{1,i}}-d_{\bar{V}_{1,i+1}}\geq d_{\bar{Z}_{1,i}}$,
which we have imposed as a condition of the theorem (and which is 
automatically satisfied if
one starts with a pair $(V_1,V_n)$ by Lemma \ref{lem-2}).

The situation for the $\bar{V}_{n,i}$ is the same,
contributing dimensions of 
$$(d_{\bar{V}_{n,i}}-d_{\bar{V}_{n,i-1}})(d_{f_{i,n-1}(E_i)}-
(d_{\tilde{Z}_i}-d_{\bar{Z}_{n,i}})-d_{\bar{V}_{n,i}})$$ 
at each step, and non-empty if and only if (\ref{2-cond-7}) is satisfied.

Finally, we need to choose $V_1$ containing $\bar{V}_{1,2}$ and
mapping surjectively onto $V_{n,1}$ under $f_{1,n-1}$. This is open
inside a Grassmannian of dimension 
$(r-d_{\bar{V}_{1,2}})(d_{V_{n,1}}+d_{\ker f_1}-r)$, with the
Grassmannian non-empty if (following our notational conventions)
$d_{\bar{V}_{1,2}} \leq d_{\bar{V}_{1,1}} \leq d_{g_{0,1}(E_1)}$, and the
open condition non-empty if
$d_{\bar{V}_{1,2}}+d_{V_{n,1}} \leq d_{\bar{V}_{1,1}}$; together, they
are non-empty if and only if (\ref{2-cond-8}) is satisfied.
Similarly, choosing $V_n$ is open inside a Grassmannian of dimension
$(r-d_{\bar{V}_{n,n-1}})(d_{V_{1,n}}+d_{\ker g_{n-1}}-r)$, non-empty
when (\ref{2-cond-9}) is satisfied. 
This describes $X$ completely, and adding up the dimension formulas and
cancelling terms gives the asserted formula and completes the proof of the 
theorem.
\end{proof} 

\begin{cor}\label{fiber-dim} The dimension of a fiber of $\pr_{1n}$ over a 
point $(V_1,V_n)$ is at most equal to 
\begin{multline*}r(d-r)
-\d(V_{1,n})(\d(\bar{V}_{1,n})-\d(V_{1,n}))
-\d(V_{n,1})(d(\bar{V}_{n,1})-d(V_{n,1})) \\
+ \sum_{i=2}^{n-1} \d(\bar{Z}_{1,i}\cap\bar{Z}_{n,i})
(\d(\tilde{Z}_i)-\d(\bar{Z}_i)) \\
- \sum_{i=1}^{n} \left(\d(\bar{Z}_{1,i})
(\d(g_{i-1,1}(E_i))-\d(\bar{V}_{1,i}))+
\d(\bar{Z}_{n,i})(\d(f_{i,n-1}(E_i))-\d(\bar{V}_{n,i}))\right). 
\end{multline*}
\end{cor}

\begin{proof} This follows immediately from the two theorems together with 
Theorem \ref{lg-main}, which implies that the dimension of $LG$ is $r(d-r)$.
\end{proof}

\section{Applications to limit linear series}\label{s-apps}

In this section, we describe the promised applications to the theory of 
limit linear series on curves. First, we recall the basic theorems on 
spaces of limit series, and how the linked Grassmannian is used to construct 
them (see \cite[\S 5]{os8} for the general case, and additional details).

We state the general theorems first, and then recall in more detail the
situation for a reducible curve over a field.

\begin{sit} Let $X/B$, together with smooth sections $P_1,\dots,P_n$, be
a {\bf smoothing family}: $B$ should be regular, and $X$ should be flat and 
proper over $B$, with fibers which are at worst nodal curves; for the
full technical details, see \cite[Def. 3.1]{os8}. We further assume that
$X/B$ has at most one node. 
\end{sit}

In \cite{os8}, for integers $r,d$ and ramification sequences $\alpha^1,
\dots,\alpha^n$, we describe a functor $\cG^r_d$ of relative limit series
of degree $d$ and dimension $r$ on $X/B$, ramified to order at least
$\alpha^i$ at each $P_i$. This functor is compatible with base change, and
agrees with usual linear series on the smooth fibers of $X/B$. We have
the following basic theorem:

\begin{thm}\label{grd-main} \cite[Thm. 5.3]{os8} In the above situation, 
the functor 
$\cG^r_d(X)$ is represented by a scheme $G^r_d(X)$, compatible with base 
change to any other smoothing family. This scheme is projective over $B$, 
and if it is 
non-empty, the local ring at any point $x \in G^r_d(X)$ closed in its fiber 
over $b \in B$ has dimension at least $\dim \cO_{B,b} + \rho$, where 
$\rho=(r+1)(d-r)-rg-\sum_{i,j}\alpha^i_j$.
\end{thm}

We now specialize to the case that that $X$ is over $\Spec k$, with two 
smooth components $Y$ and $Z$
glued at a single node $\Delta'$. We fix smooth points $P_1,\dots,P_n$ of 
$X$, as well as integers $r,d$, and 
ramification sequences $\alpha^1,\dots,\alpha^n$. In this situation, our
functor roughly parametrizes line bundles $\cL$ on $X$ of degree $d$ on $Y$
and degree $0$ on $Z$, together with vector spaces $V_i \subset H^0(X,\cL^i)$
of dimension $r$, where:
\begin{ilist}
\itm We define $\cL^i$ to be the line bundle obtained by gluing together 
$(\cL|_Y)(-i\Delta')$ and $(\cL|_Z)(i\Delta')$;
\itm Each $V_i$ maps into $V_{i+1}$ under the natural map 
$\cL^i \to \cL^{i+1}$ induced by the natural inclusion on $Z$ and the zero 
map on $Y$, and each $V_{i+1}$ maps similarly into $V_i$.
\end{ilist}

We then construct the space $G^r_d(X)$ as a closed subscheme of a 
linked Grassmannian as follows.

First, choose a divisor $D$ of very large degree, supported non-trivially on 
both $Y$ and $Z$, and disjoint from the $P_i$ and $\Delta'$. Although the 
construction depends on the choice of $D$, the resulting scheme represents
the functor $\cG^r_d(X)$, which is described independently of $D$. For 
$i=0,\dots,d$ let 
$P^i:=\Pic^{d-i,i}(X) \cong \Pic^{d-i}(Y)\times \Pic^{i}(Z)$ denote the 
Picard scheme of line bundles on $X$ restricting to degree $d-i$ on $Y$ and
degree $i$ on $Z$. Let $\cL^i$ be the universal line bundle on 
$P^i \times X$, and fix isomorphisms between the $P^i$ so that we can 
consider the $\cL^i$ as line bundles on a single scheme $P \times X$. We 
then have maps $\cL^i \to \cL^{i+1}$ and $\cL^{i+1} \to \cL^i$
as described above.

For $i=0,\dots,d$, we thus obtain maps $f_i:\cE_i \to \cE_{i+1}$ and 
$g_i:\cE_{i+1} \to \cE_i$, where $\cE_i:=p_{1*} (\cL^{i-1}(D))$, with 
$p_1:P \times X \to P$ the projection map. We then get a linked Grassmannian
$LG$ of length $n'=d+1$ over $P$, looking at sub-bundles of the $\cE_i$ of 
rank $r'=r+1$. We denote by $d'$ the rank of the $\cE_i$, so that $LG$ has
relative dimension $r'(d'-r')$ over $P$. Writing $D=D^Y+D^Z$ where $D^Y$ and
$D^Z$ are supported on $Y$ and $Z$ respectively, we then obtain our $G^r_d(X)$
space as a closed subscheme of $LG$ by requiring that the space of sections 
in $\cE_1=p_{1*}(\cL^{0}(D))$
vanishes along $D^Y$, and the space in $\cE_{n'}=p_{1*}(\cL^{d}(D))$ vanishes
along $D^Z$, which together imply that all sections will vanish along $D$,
and hence come from the original $\cL^i$. Finally, the ramification 
conditions at the $P_i$ are likewise imposed on $\cE_1$ or $\cE_{n'}$ 
depending on whether $P_i$ lies on $Y$ or $Z$. 

Given a point of $G^r_d(X)$ described by $\cL$ and $V_0,\dots,V_d$, we obtain 
a pair, which we denote by $(V^Y,V^Z)$, of linear series of degree $d$ and 
dimension $r$ on $Y$ and $Z$, by taking $V^Y=V_0|_Y$, $V^Z=V_d|_Z$ (as 
spaces of sections of the line bundles $\cL^Y:=\cL^0|_Y,\cL^Z:=\cL^d|_Z$,
which we omit from the notation). We thus have as natural map
$FR: G^r_d(X) \to G^r_d(Y) \times G^r_d(Z)$. 
The limit linear series defined by Eisenbud and Harris can be considered to
be the closed subscheme $G^{r,\EH}_d(X) \subset G^r_d(Y) \times G^r_d(Z)$
cut out by the condition that, if $a^Y_j$ and $a^Z_j$ denote the vanishing
sequences at the node of $V^Y$ and $V^Z$ respectively, for $j=0,\dots,r$
we have
\begin{equation}\label{eh-ineq}a^Y_j+a^Z_{r-j} \geq d.\end{equation}
An Eisenbud-Harris limit series is refined when these inequalities are all
equalities. In \cite{os8} the following theorem is proved:

\begin{thm}\label{grd-eh} \cite[Prop. 6.6, Cor. 6.8]{os8} The natural map 
$FR: G^r_d(X) \to G^r_d(Y) \times G^r_d(Z)$ has
set-theoretic image consisting precisely of the space of Eisenbud-Harris
limit series.
This map is an isomorphism when restricted to the open
subset of $G^r_d(X)$ mapping to refined Eisenbud-Harris limit series.
\end{thm}

In particular, we see that
$G^r_d(X)$ is obtained from $LG$ simply by imposing conditions on the first
and last projection maps, so our computations of fiber dimension for $LG$
will also apply to $G^r_d(X)$, with fibers of $\pr_{1n}$
for $LG$ corresponding precisely to fibers of $FR$ for $G^r_d(X)$.
Since we are working with such fibers, although $LG$ will be over the 
non-trivial base $P$, to study any given fiber of $\pr_{1n}$ we can 
restrict to the corresponding point on the base, which is equivalent to
fixing our choice of the line bundle $\cL$.
We now describe the relationship between
the various numerical invariants in the two situations.

\begin{lem}\label{numbers}
In the situation described above, and using the notational
conventions of Theorem \ref{crude-loci}, for a given Eisenbud-Harris
limit series $(V^Y,V^Z)$, with vanishing sequences $a^Y$ and $a^Z$ at 
$\Delta'$, we have the following dimension formulas for the corresponding
pair $(V_1,V_{n'})\in G_1 \times G_{n'}$:
$$\d(\bar{V}_{1,i})= \#\{j:a^Y_j \geq i-1\} \text{ for }i=1,\dots,n';$$
$$\d(\bar{V}_{n',i})= \#\{j:a^Z_j \geq n'-i\} \text{ for }i=1,\dots,n';$$
$$\d(g_{i-1,1}(E_i))=d+\deg D^Y+1-g_Y-i+1 \text{ for }i=2,\dots,n';$$
$$\d(f_{i,n'-1}(E_i))=d+\deg D^Z+1-g_Z-(n'-i) \text{ for }i=1,\dots,n'-1;$$
$$\d(\tilde{Z}_i)=1 \text{ for }i=2,\dots,n-1.$$
In addition, we have the following statements at the boundaries:
$$\d(V_{1,n'})= \#\{j:a^Z_j \leq 0\};$$
$$\d(V_{n',1})= \#\{j:a^Y_j \leq 0\};$$
$$\d(g_{0,1}(E_1))=d+\deg D^Y+1-g_Y-1+\d(V_{n',1});$$
$$\d(f_{n',n'-1}(E_{n'}))=d+\deg D^Z+1-g_Z-1+\d(V_{1,n'}).$$
\end{lem}

\begin{proof} We first note that $(V^Y,V^Z)$ completely determine $V_0$ and 
$V_d$ on $X$: $\cL^0|_Z$ has degree $0$, so any section in $V^Y$ vanishing 
at $\Delta'$ extends over $Z$ only by zero, while if a section is 
non-vanishing at $\Delta'$, we have $a^Y_0=0$, so $a^Z_r=d$ by 
(\ref{eh-ineq}), and we must have $\cL^Z\cong \cO_Z(d\Delta')$, so $\cL^0|_Z$
is the trivial bundle, and we can extend uniquely by a constant. The same
argument shows that $V^Z$ uniquely determines $V_d$. 

Now, by definition, for all $i$ we have that $\bar{V}_{1,i}$ is the space of 
sections of $V^0$ in the image of $H^0(X,\cL^{i-1}(D))$, which is exactly
the space of sections of $\cL$ vanishing to order at least $i-1$ at 
$\Delta'$, giving the asserted dimension. We defined 
$V_{1,n'}$ to be the space of sections of $V^0$ in the image of $V^d$;
we note that this is zero-dimensional if $V^d$ has no sections non-vanishing
on $Y$, and one-dimensional otherwise, giving the asserted formula
for $\d(V_{1,n'})$. The formulas for $\d(\bar{V}_{n',i})$ and
$\d(V_{n',1})$ are obtained similarly.

Next, $g_{i-1,1}(E_i)$ for $i\geq 2$ are those sections of $E_1$ in the image
of $E_i$; such sections necessarily vanish on $Z$, and are in fact in 
correspondence with sections of $\cL^0(D)|_Y$ vanishing to order at least
$i-1$ at $\Delta'$. Because $D$ was chosen sufficiently large, 
$h^0(Y,\cL^0(D)|_Y)=\deg \cL^0(D)|_Y+1-g_Y=d+\deg D^Y+1-g_Y$, and the number
of sections vanishing to order at least $i-1$ at $\Delta'$ is
$d+\deg D^Y+1-g_Y-i+1$, as asserted. For $i=1$, by our notational conventions
we need to consider $\d(\bar{V}_{n,1})+\d (\ker f_1)$. We have 
$\ker f_1=g_1(E_1)$, so we know this has dimension $d+\deg D^Y+1-g_Y-1$, and
we obtain the asserted formula.
The same argument works to compute $\d(f_{i,n'}(E_i))$.  

Finally, we have $\d(\tilde{Z}_i) \leq 1$ for all $i$ because given a line 
bundle $\cL$ on $X$, the only obstruction to gluing sections of $\cL|_Y$ and 
$\cL|_Z$ to obtain a section of $\cL$ is whether the sections agree in 
$\cL|_{\Delta'} \cong k$. On the other hand, we have 
$\d(\tilde{Z}_i) \geq 1$ because 
for $\cL=\cL^i(D)$, both restriction maps $\cL|_Y \to \cL_{\Delta'}$ and
$\cL|_Z \to \cL_{\Delta'}$ are surjective because $D$ was chosen to be 
sufficiently ample.
\end{proof}

\begin{cor}\label{crude-dim} The dimension of the space of crude limit 
series corresponding
to a given Eisenbud-Harris limit series $(V^Y,V^Z)$ is bounded above by
$$\sum_{i=0}^r(a^Y_i+a^Z_{r-i}-d).$$
\end{cor}

\begin{proof}
This is a direct application of Corollary \ref{fiber-dim},
together with the preceding lemma. We begin by recasting the formula
\begin{multline*}\sum_{i=1}^{n} \left(\d(\bar{Z}_{1,i})
(\d(g_{i-1,1}(E_i))-\d(\bar{V}_{1,i}))+
\d(\bar{Z}_{n',i})
(\d(f_{i,n'-1}(E_i))-\d(\bar{V}_{n',i}))\right) \\
+\d(V_{1,n'})(\d(\bar{V}_{1,n'})-\d(V_{1,n'}))
+\d(V_{n',1})(d(\bar{V}_{n',1})-d(V_{n',1})) \\
- \sum_{i=2}^{n'-1} \d(\bar{Z}_{1,i}\cap\bar{Z}_{n',i})
(\d(\tilde{Z}_i)-\d(\bar{Z}_i)),
\end{multline*}
in terms of the numerical invariants of $(V^Y,V^Z)$. 

We first see from the lemma that $\d(\tilde{Z}_i)=1$, and since 
$\d(\bar{Z}_{1,i}\cap\bar{Z}_{n',i}) \leq \d(\bar{Z}_i)$, the term
$\sum_{i=2}^{n'-1} \d(\bar{Z}_{1,i}\cap\bar{Z}_{n',i})
(\d(\tilde{Z}_i)-\d(\bar{Z}_i))$ always vanishes. Similarly, we have
$\d(\bar{V}_{1,n'})=\#\{j:a^Y_j \geq d\}$, so is at most $1$, and likewise
$\d(\bar{V}_{n',1}) \leq 1$, so we have that the terms
$\d(V_{1,n'})(\d(\bar{V}_{1,n'})-\d(V_{1,n'}))
+\d(V_{n',1})(d(\bar{V}_{n',1})-d(V_{n',1}))$
always vanish. Thus, it is enough to consider the first sum.

We now claim that we have:
\begin{multline*} \sum_{i=1}^{n} \d(\bar{Z}_{1,i}) 
(\d(g_{i-1,1}(E_i))-\d(\bar{V}_{1,i}))+
\d(\bar{Z}_{n,i})
(\d(f_{i,n-1}(E_i))-\d(\bar{V}_{n,i})) \\
= r'(d'-r') - \sum_{j=0}^r(a^Y_j+a^Z_{r-j}-d).
\end{multline*}

Indeed, we see that
$$\d(\bar{Z}_{1,i})=\d(\bar{V}_{1,i}) - \d(\bar{V}_{1,i+1}) =
\begin{cases} 1: \exists j \text{ such that } a^Y_j=i-1 \\
0: \text{otherwise}
\end{cases},$$
so if we split the sum in two, the first sum may 
be rewritten as
$$\delta_Y+\sum_{j=0}^r (d+\deg D^Y+1-g_Y-a^Y_j-(r+1-j)),$$
and similarly the second sum is
$$\delta_Z+\sum_{j=0}^r (d+\deg D^Z+1-g_Z-a^Z_j-(r+1-j)),$$
where the $\delta_Y$ and $\delta_Z$ account for the possible discrepency
between 
$\d(g_{0,1}(E_1))$ and $d+\deg D^Y+1-g_Y$ and between
$\d(f_{n',n'-1}(E_{n'}))$ and $d+\deg D^Z+1-g_Z$ respectively. However,
we claim that in fact $\delta_Y=\delta_Z=0$. Indeed, applying the boundary
cases of the lemma we have that
$\d(g_{0,1}(E_1))=d+\deg D^Y+1-g_Y$ unless $\d(V_{n,1})=0$, in which
case $a^Y_0>0$. But the term 
$(\d(\bar{V}_{1,1}) - \d(\bar{V}_{1,2})) 
(\d(g_{0,1}(E_1))-\d(\bar{V}_{1,1}))$ only contributes to the sum if some
$a^Y_j=0$, or equivalently, if $a^Y_0=0$. Thus, this term only appears in
the sum if $\d(g_{0,1}(E_1))=d+\deg D^Y+1-g_Y$, and we have $\delta^Y=0$.
The same argument shows that $\delta^Z=0$.

Next, noting that we can replace $(r+1-j)$ in the first sum by $(1+j)$, 
combining them we get
\begin{multline*}
(r+1)(d+\deg D^Y+\deg D^Z -g_Y-g_Z)-\sum_{j=0}^r(a^Y_j+a^Z_{r-j}-d+2j) \\
= r'(d+\deg D -g)-r(r+1)-\sum _{j=0}^r(a^Y_j+a^Z_{r-j}-d) \\
= r'(d+\deg D + 1 - g - (r+1))-\sum_{j=0}^r(a^Y_j+a^Z_{r-j}-d) \\
= r'(d'-r') - \sum_{j=0}^r(a^Y_j+a^Z_{r-j}-d),
\end{multline*}
as desired.

Finally, putting these statements together and applying Corollary 
\ref{fiber-dim}, we obtain the desired bound.
\end{proof}

This allows us to conclude that if the components of a reducible curve have 
the expected dimensions of (ramified) linear series, then the curve as a 
whole has the expected dimension of limit linear series:

\begin{cor}\label{lim-cor} Let $X$ be a proper curve of genus $g$ over 
$\Spec k$, consisting
of the union of two smooth components $Y$ and $Z$ at a point $\Delta'$. Also 
fix distinct marked points $P_1,\dots,P_n$ in the smooth locus of $X$. Given
integers $d,r$, and ramification sequences 
$\alpha^i:=\alpha^i_0,\dots,\alpha^i_r$, suppose that the following 
condition holds:

For any $\alpha^{\Delta'}:=\alpha^{\Delta'}_0,\dots,\alpha^{\Delta'}_r$, the 
space of linear series
on $Y$ of degree $d$ and dimension $r$ with ramification sequence at least
$\alpha^i$ at any $P_i$ lying on $Y$, and at least $\alpha^{\Delta'}$ at 
${\Delta'}$, has the expected dimension $\rho^Y:=(r+1)(d-r)-rg_Y-
\sum_{P_i\in Y,j} \alpha^i_j-\sum_j \alpha^{\Delta'}_j$ if it is non-empty, 
and similarly for $Z$.

Then the space $G^r_d(X)$ of limit linear series on $X$ is pure of dimension 
$\rho^X:=(r+1)(d-r)-rg-\sum_{i,j} \alpha^i_j$, and non-empty if and only
if there exist ramification sequences $\alpha^Y$ and $\alpha^Z$ with
$\alpha^Y_j+\alpha^Z_{r-j}=d-r$ for each $j$ such that the corresponding
spaces $G^r_d(Y)$ and $G^r_d(Z)$ are both non-empty.
\end{cor}

\begin{proof} By Theorem \ref{grd-eh}, we know that the set-theoretic 
image of $G^r_d(X)$ in the space 
$G^r_d(Y) \times G^r_d(Z)$ is precisely the closed subscheme of pairs
whose vanishing sequences $a^{\Delta',Y}_j$ and $a^{\Delta',Z}_j$ at $\Delta'$ satisfy the
inequalities
$$a^{\Delta',Y}_j + a^{\Delta',Z}_{r-j} \geq d,$$
or equivalently, whose ramification sequences satisfy
$$\alpha^{\Delta',Y}_j + \alpha^{\Delta',Z}_{r-j} \geq d-r.$$
If we fix ramification sequences meeting this condition, we find that the
dimension of the locus in $G^r_d(Y)\times G^r_d(Z)$ is at most 
\begin{multline*}\rho^Y+\rho^Z =  \\
(r+1)(d-r)-rg_Y- \sum_{P_i\in Y,j} \alpha^i_j-\sum_j \alpha^{\Delta',Y}_j \\
+ (r+1)(d-r)-rg_Z- \sum_{P_i\in Z,j} \alpha^i_j-\sum_j \alpha^{\Delta',Z}_j \\
= 2(r+1)(d-r)-r(g_Y+g_Z) -\sum_{i,j} \alpha^i_j -\sum_j (\alpha^{\Delta',Y}_j+
\alpha^{\Delta',Z}_{r-j}) \\
= \rho^X -\sum_j(\alpha^{\Delta',Y}_j+\alpha^{\Delta',Z}_{r-j}-(d-r)) \leq 
\rho^X
\end{multline*}
because of the above inequality and the identity $g_Y+g_Z=g$. 
We see further that to show that $\dim G^r_d(X) \leq \rho^X$, it suffices to
see that the dimension of the fibers of
the map $G^r_d(X) \rightarrow G^r_d(Y)\times G^r_d(Z)$ are at most
$\sum_j(\alpha^{\Delta',Y}_j+\alpha^{\Delta',Z}_{r-j}-(d-r)
)=\sum_j(a^Y_j+a^Z_{r-j}-d)$,
which is precisely Corollary \ref{crude-dim}.
But we have $\dim G^r_d(X) \geq \rho^X$ by Theorem \ref{grd-main},
so we obtain equality. Finally, for the non-emptyness assertion, 
we note that even if the only Eisenbud-Harris limit series on $X$ are crude,
we can impose weaker ramification conditions at $\Delta'$ to satisfy the desired
equality, and we will still have that the corresponding spaces $G^r_d(Y)$ 
and $G^r_d(Z)$ are non-empty.
\end{proof}

The base cases for our induction will be the following easy and
well-known statements:

\begin{lem}\label{base-case} Let $C$ be a smooth proper curve of genus $g$ over $\Spec k$, 
with distinct marked points $P_1,\dots,P_n$. Then for any $r,d$, and
ramification sequences $\alpha^i:=\alpha^i_0,\dots,\alpha^i_r$, the
space of linear series of degree $d$ and dimension $r$ on $C$ having
ramification sequences at least $\alpha^i$ at $P_i$ for all $i$ has
dimension exactly $\rho:=(r+1)(d-r)-rg-\sum_{i,j} \alpha^i_j$ (but is
not necessarily non-empty) in either of the following cases:
\begin{ilist}
\itm if $g=0$ and $\ch k=0$;
\itm if $g=1$ and $n=1$.
\end{ilist}
Assuming $\rho \geq 0$, case (ii) is non-empty if and only if the single 
imposed vanishing sequence does not have $a_r=d$, $a_{r-1}=d-1$. 
Non-emptiness of case (i) is determined by Schubert calculus, but in 
particular is non-empty whenever every imposed ramification sequence is 
of the form $0,1,\dots,1$.
\end{lem}

Note that no generality hypotheses are required here. (i) fails in positive
characteristic even for $r=1$, although it may be proved in certain 
non-trivial cases.

\begin{proof} The dimension statement for (i) is proved by an inductive 
argument and the Pl\"ucker formula; see \cite[Thm. 2.3]{e-h4}. Non-emptiness
is simply a question of Schubert calculus, and we recall the argument: 
the ramification 
conditions of the form $0,1,\dots,1$ correspond to an intersection of 
Schubert classes $\sigma_{0,1,\dots,1}$ in the Grassmannian $\GG(r,d)$.
If we pass to the dual, we obtain a collection of special Schubert classes
$\sigma_{0,\dots,0,r}$ in $\GG(d-r-1,d)$, and can easily check non-emptiness
as long as the expected dimension is non-negative by inductively
applying Pieri's formula \cite[p. 271]{fu1}.

(ii) may be seen as follows: 
let the prescribed vanishing sequence at $P_1$ be $a_0,\dots,a_r$. The
expected dimension is $\rho=(r+1)(d-r)-r-\sum_{i=0}^r(a_i-i)$. 

We first consider the case that $a_r=d$. Then the only possibility is a
linear series contained in $H^0(C,\cO(dP_1))$. This space is $d$-dimensional,
with sections vanishing to every order at $P_1$ except $d-1$. Thus the
space of linear series is contained in a Grassmannian $\GG(r,d-1)$, of 
dimension $(r+1)(d-r-1)$, and the ramification condition cuts out a 
Schubert cycle of codimension $\sum_{i=0}^{r-1}(a_i-i)+(a_r-r-1)$ as long as 
$a_{r-1}<d-1$; if $a_{r-1}=d-1$, the space is necessarily empty. We thus
obtain the desired statement in this case.

If $a_r<d$, the line bundle could be any line bundle of degree $d$; these
are all of the form $\cL =\cO((d+1)P_1-Q)$ as $Q$ varies over the points of 
$C$. We have $h^0(C,\cL)=d$, and there are sections
vanishing to all orders less than $d-1$ at $P_1$. When $Q=P_1$, the last 
order of 
vanishing is $d$, while for $Q \neq P_1$, the last is $d-1$. The ramification
condition is imposed inside a $\GG(r,d)$-bundle over $\Pic(C) \cong C$, of
dimension $1+(r+1)(d-r-1)$; since we only need an upper bound on the 
dimension, we may work fiber by fiber. We claim that for each $\cL$, 
the ramification condition imposes a Schubert cycle of
codimension $\sum_{i=0}^r(a_i-i)$, giving the correct dimension for the 
total space. This is clear if $a_r<d-1$, or if $a_r=d-1$ and $Q \neq P_1$. 
In the last case, we note that the condition imposed by $a_r=d-1$ is the
same as that imposed by $a_r=d$, so we still get a Schubert cycle of the
asserted dimension. 

Thus, we see that the space always has the asserted dimension, and is
non-empty as long as we do not have $a_r=d,a_{r-1}=d-1$.
\end{proof}

We are now ready for:

\begin{proof}[Proof of Theorem \ref{main}] As promised, everything except 
the connectedness statement will 
follow immediately from the limit linear series machinery; we will conclude
connectedness in the reducible case from the theorem of Fulton and Lazarsfeld
in the irreducible case. We first note that properness and the
lower bound on the dimension follow directly from Theorem \ref{grd-main},
and require no generality or characteristic hypotheses. We 
will prove the upper bound on dimension, and non-emptiness by induction
on $g$. 

In fact, we induct on a slightly stronger statement: we will show 
non-emptiness for $\rho \geq 0$ also in the case that we have imposed
ramification sequences $0,1,\dots,1$ at general points.
The base case is $g=0$; this is case (i) of Lemma \ref{base-case}.
Suppose we now know the statement for curves of genus less than $g$, and we 
want to conclude it for curves of genus $g$. We first note that the
reducible case for any curve with both components of genus strictly less than
$g$ follows immediately, by Corollary \ref{lim-cor}, and the reducible
case in full generality will likewise follow once we have proved the 
irreducible case for genus $g$. To prove this case, we consider specifically
a curve $X_0$ consisting of one component $Y$ having genus $1$ and no imposed 
ramification points, and the other component $Z$ of genus $g-1$, with all 
imposed ramification points. By case (ii) of Lemma \ref{base-case} and by our
induction hypothesis, we have the dimensional upper bound on each component
for any ramification sequences at the node, and non-emptiness when we 
consider the vanishing sequences $d-r-1,d-r,\dots,d-2,d$ on $Y$, and 
$0,2,3,\dots,r+1$ on $Z$, so we obtain both statements for $G^r_d(X_0)$
by Corollary \ref{lim-cor}. We place $X_0$ in a smoothing family $X/B$ 
\cite[Thm. 3.4]{os8} with smooth generic fiber, and because the 
corresponding space
$G^r_d(X)$ of relative limit series is proper, we conclude the dimensional
upper bound and non-emptiness statements for the generic fiber, which is
enough to imply them for a general curve, since $\cM_{g,n}$ is connected.

Finally, in the case $\rho>0$, we show connectedness. Fulton and Lazarsfeld
\cite{f-l1} proved 
connectedness in the irreducible case. If we start with a reducible curve
$X_0$ satisfying the hypotheses of our theorem, we place it as before in a 
smoothing family $X/B$ with smooth generic fiber, and regular one-dimensional 
base. By Fulton-Lazarsfeld, the space of linear series is connected over the
open subset of the base corresponding to smooth curves. Because $X_0$ is 
general, $G^r_d(X_0)$ has dimension $\rho$, and no component of $G^r_d(X)$ is 
supported over the special fiber. Because $G^r_d(X)$ is proper, we conclude 
also connectedness of the special fiber by, e.g., \cite[Prop. 15.5.3]{ega43}.
\end{proof}

\begin{rem} In fact, we can use the same sort of arguments to give a
sharp statement in terms of Schubert calculus on when $G^r_d$ spaces
are non-empty. However, this statement is already known \cite[Rem. following
Thm. 4.5]{e-h1}, and we do not pursue it.
\end{rem}

We conclude with a simple lemma in the limit linear series context. This
lemma serves both to show that the upper bound of Corollary \ref{crude-dim}
is not sharp, and also to show that in the case of
``codimension 1'' crude limit series, we still have a set-theoretic 
equivalence between limit series and Eisenbud-Harris limit series. 

\begin{lem} Continuing with the notation of Corollary \ref{lim-cor},
let $(V^Y,V^Z)$ be a pair of linear series on $Y$ and $Z$, each of 
degree $d$ and dimension $r$, satisfying the desired ramification conditions,
and such that for some $j_0$, we have
$$d+1 \geq a^{\Delta',Y}_{j_0} + a^{\Delta',Z}_{r-j_0} \geq d,$$
while for other $j$ we still have 
$a^{\Delta',Y}_j + a^{\Delta',Z}_{r-j} = d$.
Then there is a unique limit linear series on $X$ restricting to $(V^Y,V^Z)$.
\end{lem}

\begin{proof} In the case that 
$$a^{\Delta',Y}_j + a^{\Delta',Z}_{r-j} = d$$
for all $j$,
the arguments of \cite[\S 6]{os8} come out of the following observation:
this equality means that for any $i$, if we denote by $d_i^Y$ the dimension 
of the space of sections of $V^Y$ vanishing to order at least $i$ at $\Delta'$,
and $d_i^Z$ the same for $Z$, then we have 
$$d_i^Y+d_i^Z= \begin{cases} r+1:\nexists j \text{ such that } 
i=a^{\Delta',Y}_j=d-a^{\Delta',Z}_{r-j}\\ 
r+2: \text{ otherwise.} \end{cases}$$
For $i$ such that $i \neq a^{\Delta',Y}_j$ for any $j$, all these sections vanish
at $\Delta'$ on both $Y$ and $Z$, and we can glue them arbitrarily in $\cL^i$, so we 
get an $(r+1)$-dimensional space of possible sections of $\cL^i$, and we are 
forced to choose this for $V_i$. If $i= a^{\Delta',Y}_j$ for some $j$, in principal 
we are choosing a $V_i$ contained in an $(r+2)$-dimensional space, but here
we have sections which are non-vanishing at $\Delta'$ on both $Y$ and $Z$, so the 
requirement that they glue together at $\Delta'$ is a codimension $1$ condition, 
and once again we are choosing $V_i$ from an $(r+1)$-dimensional space. 

But now suppose that we also allow
$$a^{\Delta',Y}_{j_0} + a^{\Delta',Z}_{r-j_0} = d+1$$
for a certain $j_0$. 
The basic observation is that for the above argument to work, it suffices
to know that $d_i^Y+d_i^Z \leq r+2$, with equality only when the appropriate
space of sections on either $Y$ or $Z$ has sections non-vanishing at $\Delta'$,
since we still obtain a codimension 1 gluing condition in this case. But
this is equivalent to having $d_i^Y+d_i^Z=r+2$ only when $i=a^{\Delta',Y}_j$
or $i=d-a^{\Delta',Z}_{r-j}$ for some $j$, which one easily sees will follow from
our hypothesis.
\end{proof}

We can thus conclude:

\begin{cor} Continuing with the hypotheses and notation of 
Theorem \ref{main}, if $\rho \leq 1$ we have that the natural set-theoretic
map $G^r_d(X) \to G^{r,\EH}_d(X)$ is a bijection.
\end{cor}

\begin{proof} We know by Theorem \ref{grd-eh} that the map is always 
surjective, so it is enough to 
prove injectivity. Let $(V^Y,V^Z)$ be a point of $G^{r,\EH}_d(X)$, with
vanishing sequences $a^Y_i$ and $a^Z_i$ at the node. 
Because everything is general, for $(V^Y,V^Z)$ to exist, by Theorem 
\ref{main} we must have that
$\rho_Y$ and $\rho_Z$ are both non-negative, after taking into account the
ramification at the node. Furthermore, by the additivity of the 
Brill-Noether number we have that 
$$\rho-\rho_Y-\rho_Z=\sum_j(a^Y_j+a^Z_{r-j}-d).$$
Since $\rho \leq 1$ and $\rho_Y$, $\rho_Z$, and the right hand side are
all at least $0$, we see that the right-hand side is either $0$ or $1$.
The case that it is $0$ is the case that $(V^Y,V^Z)$ is refined, in which
case we already knew that there is a unique point of $G^r_d(X)$ above it,
and the case that it is $1$ is the case addressed by the lemma.
\end{proof}

\bibliographystyle{hamsplain}
\bibliography{hgen}

\end{document}